\documentclass{amsart}

\def\C{\mathbb{C}}

\def\l{\lambda}
\def\bl{\boldsymbol}

\def\D{\mathbb D}

\def\ov{\overline}
\def\lo{\longrightarrow}

\def\mb{\mathbb}

\def\d{\displaystyle\sum}

\def\i{\prime}

\newcommand{\be}{\begin{equation}}
\newcommand{\ee}{\end{equation}}
\newcommand{\bea}{\begin{eqnarray}}
\newcommand{\eea}{\end{eqnarray}}
\newcommand{\Bea}{\begin{eqnarray*}}
\newcommand{\Eea}{\end{eqnarray*}}
\theoremstyle{definition}

\theoremstyle{remark}

\numberwithin{equation}{section}

\newtheorem{thm}{Theorem}[section]
\newtheorem{cor}[thm]{Corollary}
\newtheorem{lem}[thm]{Lemma}
\newtheorem{prop}[thm]{Proposition}

\begin{document}

\title[Kernel functions on the symmetrized polydisc]{Reproducing kernel for
a class of weighted Bergman spaces on the symmetrized polydisc}



\author[Misra]{Gadadhar Misra}
\address[Gadadhar Misra]{Department of Mathematics, Indian Institute of Science,
Bangalore 560 012, India}
\email[Misra]{gm@math.iisc.ernet.in}

\author[Shyam Roy]{Subrata Shyam Roy}
\address[Subrata Shyam Roy]{Indian Institute of Science Education and Research,
Kolkata, Mohanpur Campus, Mohanpur (West Bengal) - 741 252, India}
\email[Subrata Shyam Roy]{ssroy@iiserkol.ac.in}

\author[Zhang]{Genkai Zhang}
\address[Genkai Zhang]{Department of Mathematics\\
Chalmers University of Technology \vskip 0.0ex
and Gothenburg University, S-412 96 Gothenburg,
Sweden}
\email[Genkai Zhang]{genkai@chalmers.se}

\thanks{Financial support for the work of G. Misra and Genkai Zhang was provided
by the Swedish Research Links programme entitled ``Hilbert modules, operator
theory and complex analysis".}
\date{}

\dedicatory{}

\commby{}

\begin{abstract}
A natural class of weighted Bergman spaces on the symmetrized polydisc is
isometrically embedded as
a subspace in the corresponding weighted Bergman space on the polydisc.  We find
an orthonormal basis
for this subspace. It enables us to  compute the kernel function for the
weighted Bergman spaces on the
symmetrized polydisc using the explicit nature of our embedding. This family of
kernel functions include the
Szeg\"{o} and the Bergman kernel on the symmetrized polydisc.

\end{abstract}

\maketitle

\section{Introduction}
Let $\varphi_i,i\geq 0,$ be the \emph{elementary symmetric function} of degree
$i$, that is, $\varphi_i$ is the sum of all products of $i$ distinct variables
$z_i$ so that $\varphi_0=1$ and
$$
\varphi_i(z_1,\ldots ,z_n) = \sum_{1\leq k_1< k_2 <\ldots <k_i\leq n} z_{k_1}
\cdots z_{k_i}.
$$
For $n \geq 1,$ let $\mathbf s:\C^n\lo\C^n$ be the function of symmetrization
given by the formula
$$\mathbf s(z_1,\ldots,z_n)= \big (\varphi_1(z_1,\ldots , z_n), \ldots ,
\varphi_n(z_1,\ldots ,z_n) \big ).$$
The image $\mathbb G_n:=\mathbf s(\D^n)$ under the map $\mathbf s$ of the unit
polydisc $\mathbb D^n:=\{\bl z\in \mathbb C^n: \|\bl z\|_\infty < 1\}$ is known
as the {\it symmetrized polydisc}. The restriction map $\mathbf s_{|{\rm
res}\,\mathbb D^n}:\mathbb D^n \to \mathbb G_n$ is a proper holomorphic map
\cite{R}.
The Bergman kernel for the symmetrized polydisc is computed explicitly in
\cite{EZ}. It is obtained from the transformation rule for the Bergman kernel
under proper holomorphic maps \cite[Theorem 1]{Bell}.

Here we realize  (isometrically) the Bergman space $\mathbb A^2(\mathbb G_n)$
of the symmetrized polydisc as a subspace of the Bergman space $\mathbb
A^2(\mathbb D^n)$ on the polydisc using the symmetrization map $\mathbf s$.
Indeed, the map
$\Gamma: \mathbb A^2(\mathbb G_n)  \to \mathbb A^2(\mathbb D^n)$ defined by the
formula
$$
(\Gamma f)(\boldsymbol z) = (f\circ \mathbf s)(\boldsymbol z) J_\mathbf
s(\boldsymbol z),\,\, \boldsymbol z\in \mathbb D^n,
$$
where $J_\mathbf s$ is the complex Jacobian of the map $\mathbf s$, is an
isometric embedding.  The image ${\rm ran}\, \Gamma \subseteq \mathbb
A^2(\mathbb D^n)$ consists of anti-symmetric functions:
$$
\mbox{\rm ran}\, \Gamma:=\{f: f(\boldsymbol z_\sigma) = {\rm sgn}(\sigma)
f(\boldsymbol z),\, \sigma \in \Sigma_n \,,f\in \mathbb A^2(\mathbb D^n)\},
$$
where $\Sigma_n$ is the symmetric group on $n$ symbols. The range of $\Gamma$ is
a subspace of $\mathbb A^2(\mathbb D^n)$, we let $\mathbb A_{\rm anti}^2(\mathbb
D^n)$ be this subspace. An orthonormal basis of $\mathbb A_{\rm anti}^2(\mathbb
D^n)$ may then be transformed in to an orthonormal basis of the $\mathbb
A^2(\mathbb G_n)$
via the unitary map $\Gamma^*$. It is then possible to compute the Bergman
kernel for the symmetrized polydisc $\mathbb G_n$ by evaluating the sum
$$
\sum_{k \geq 0} e_k(\boldsymbol z)\overline{e_k(\boldsymbol w)},\,\,\boldsymbol
z,\boldsymbol w \in \mathbb G_n,
$$
for some choice of an orthonormal basis in $\mathbb A^2(\mathbb G_n)$.

This scheme works equally well for a class of weighted Bergman spaces $\mathbb
A^{(\lambda)}(\mathbb D^n)$, $\lambda >  1,$ determined by the kernel function
$$
\mathbf B_{\mathbb D^n}^{(\lambda)}(\boldsymbol z, \boldsymbol w)=
\prod_{i=1}^n(1-z_i\bar{w}_i)^{-\lambda},\,\boldsymbol z =(z_1,\ldots ,z_n),\,
\boldsymbol w=(w_1,\ldots , w_n) \in \mathbb D^n,
$$
defined on the polydisc and the corresponding weighted Bergman spaces $\mathbb
A^{(\lambda)}(\mathbb G^n)$
on the symmetrized polydisc.

The limiting case of $\lambda=1$, as is well-known, is the Hardy space on the
polydisc. We show that the reproducing kernel for the Hardy space of the
symmetrized polydisc is of the form
$$
\mathbb S_{\mathbb G_n}^{(1)}(\mathbf s(\bl z), \mathbf s(\bl w)) =
\prod_{i,j=1}^n(1-z_i\bar{w}_j)^{-1}, \bl z,\bl w\in \mathbb D^n.
$$
This is a consequence of the determinantal identity \cite[(4.3),  pp. 63]{IGM}.
Indeed, along the way, we obtain a generalization of this well-known identity.
We also point out that the Hardy kernel is not a power of the Bergman kernel
unlike the case of bounded symmetric domains.

\section{Weighted Bergman spaces on the symmetrized polydisc}
For $\lambda >1$, let $dV^{(\lambda)}$ be the probability  
measure $\big(\frac{\lambda-1}{\pi}\big )^n \Big
(\prod_{i=1}^n(1-r_i)^{\lambda-2}r_i dr_i d\theta_i\Big )$
on the polydisc $\mathbb D^n$.  Let
$dV^{(\lambda)}_\mathbf s$ be the measure on
the symmetrized polydisc $\mathbb G_n$
obtained by the change of variable formula:
$$
\int_{\mathbb G_n} f \,dV^{(\lambda)}_\mathbf s
=\int_{\D^n} (f\circ\mathbf s)\, |J_\mathbf s|^2 dV^{(\lambda)},\,\, \lambda > 1
$$
where $J_\mathbf s(\bl z) = \prod_{1\leq i < j\leq n} (z_i-z_j)$ is the complex
Jacobian of the symmetrization map
$\mathbf s$. Let  $\|J_\mathbf s\|^2_\lambda = \int_{\mathbb D^n} |J_\mathbf
s|^2 dV^{(\lambda)}$ be the the norm of the jacobian determinant $J_\mathbf s$
in the Hilbert space $L^2(\mathbb D^n, dV^{(\lambda)})$.
By a slight abuse of notation, we let $dV^{(\lambda)}_\mathbf s$ be the  measure
$\|J_s\|_\lambda^{-2}{dV^{(\lambda)}_\mathbf s}$,  $\lambda > 1$, on the
symmetrized polydisc $\mathbb G_n$.
The weighted Bergman space
$\mathbb A^{(\lambda)}(\mathbb G_n)$, $\lambda >1$, on the symmetrized polydisc
$\mathbb G_n$ is the subspace of the Hilbert space
$L^2(\mathbb  G_n,\,dV^{(\lambda)}_\mathbf s)$ consisting of holomorphic
functions.  It coincides with the usual Bergman space for $\lambda =2$.
The norm of $f\in \mb A^{(\l)}(\mathbb G_n)$ is given by
$\|f\|^2=\int_{\mathbb G_n}|f|^2 dV^{(\lambda)}_{\mathbf s}.$
We have normalized the volume measure on $\mathbb G_n$ to ensure  $\|1\|=1$.

For $\lambda >1$, let $\Gamma: \mb A^{(\l)}(\mathbb G_n)\lo \mb A^{(\l)}(\D^n)$
be the operator defined by the rule:
$$(\Gamma f)(\bl z)=  \|J_\mathbf s\|_\lambda^{-1}  J_\mathbf s(\bl
z)(f\circ\mathbf s) (\bl z),\,\,\,\,\,\,f\in \mb A^{(\l)}(\mathbb G_n),\,\, \bl
z\in\D^n.$$
It is clear from the definition of the norm in $\mb
A^{(\l)}(\mathbb G_n)$ that $\Gamma$ is an isometry. The image of $\mb
A^{(\l)}(\mathbb G_n)$ under the isometry $\Gamma$ in
$\mb A^{(\l)}(\D^n)$ is the subspace $\mb A_{\rm anti}^{(\l)}(\D^n)$ of
anti-symmetric functions since
$J_\mathbf s(\bl z_\sigma) = {\rm sgn}( \sigma) J_\mathbf s(\bl z)$, $\sigma\in
\Sigma_n$.
Every function $g$ in $\mathbb A_{\rm anti}^{(\lambda)}(\mathbb D^n)$ is of the
form
$J_\mathbf s h$ for some symmetric function $h$. For instance, take $h=J_\mathbf
s^{-1}
g$ on the open set $\{(z_1,\ldots , z_n) \in \mathbb D^n: z_i \not = z_j,\,
i\not = j\}$.
It follows that $g= J_\mathbf s (f\circ \mathbf s)$ for some function $f$
defined on
$\mathbb G_n$. Therefore, the range of the isometry coincides with the subspace
$\mathbb A_{\rm anti}^{(\lambda)}(\mathbb D^n)$. Now, it is easily verified that
$\Gamma^* g = \|J_\mathbf s\|_\lambda\,f$, where $f$ is chosen satisfying $g(\bl
z)= J_\mathbf s(\bl z) (f\circ\mathbf s) (\bl z)$.
The operator  $\Gamma:\mb A^{(\l)}(\mathbb G_n)\lo \mb A_{\rm
anti}^{(\l)}(\D^n)$ is evidently unitary.
The Hilbert spaces $\mb A^{(\l)}(\mathbb G_n)$, $\lambda > 1$, are the
weigheted Bergman spaces on the symmetrized polydisc $\mb G_n$.

Since the subspace $\mb A_{\rm anti}^{(\l)}(\D^n)$ is invariant under the
multiplication by the elementary symmetric function $\varphi_i$, $1\leq i \leq
n$, we see that it
admits a module action via the map
$$
(p,f) \mapsto p(\varphi_1, \ldots , \varphi_n) f,\,\, f \in  \mb A_{\rm
anti}^{(\l)}(\D^n),\,\, p \in \mb C[\bl z]
$$
over the polynomial ring $\mb C[\bl z]$.  The polynomial ring acts naturally via
multiplication by the coordinate functions
on the Hilbert space $\mb A^{(\l)}(\mathbb G_n)$ making it a module over the
polynomial ring $\mb C[\bl z]$.
The unitary operator $\Gamma$  intertwines the multiplication by the elementary
symmetric functions on the
Hilbert space $\mb A_{\rm anti}^{(\l)}(\D^n)$ with the multiplication
by the co-ordinate functions on $\mb A^{(\l)}(\mathbb G_n)$.  Thus
$\mb A^{(\l)}(\mathbb G_n)$ and $\mb A_{\rm anti}^{(\l)}(\D^n)$ are isomorphic
as modules via the unitary map $\Gamma$. Moreover,
since  $\mb A^{(\l)}(\mathbb G_n)$ is a submodule of the $L^2(\mathbb
G_n,dV^{(\lambda)}_{\mathbf s})$, it follows that the map
$$
(p,f) \mapsto p\cdot f,\,\, f \in  \mb A^{(\l)}(\mathbb G^n),\,\, p \in \mb
C[\bl z]
$$
is contractive. It therefore extends to a continuous map of the function algebra
$\mathcal A(\mathbb G_n)$ obtained by taking the closure
of the polynomial ring with respect to the supremum norm on the symmetrized
poly-disc.

\subsection{Orthonormal basis and kernel function}
A \emph{partition} {\boldmath${p}$} is any finite sequence {$\boldsymbol p$}
$:=(p_1,\ldots , p_n)$
of non-negative integers in decreasing order, that is,
$$
p_1 \geq \cdots \geq p_n.
$$
We let $[n]$ denote the set of all partitions of size $n$.
If a partition $\bl p$ also has the the property $p_1 >p_2 > \cdots > p_n \geq
0$, then we may write $\bl p = \boldsymbol{m} + \boldsymbol{\delta}$,
where $\boldsymbol{m}$ is some partition in $[n]$ and
$\boldsymbol{\delta}=(n-1,n-2,\ldots , 1,0)$.   Let $[\![n]\!]$ be the set of
all partitions of
the form $\bl m + \bl \delta$ for $\bl m \in [n]$.

Let ${\bl z}^{\boldsymbol m}:=z_1^{m_1}\cdots z_n^{m_n}$,
$\boldsymbol{m}\in[n]$, be a monomial.
Consider the polynomial $a_{\boldsymbol m}$ obtained by anti-symmetrizing the
monomial $\bl z^{\boldsymbol m}$:
$$
a_{\boldsymbol m}(\bl z):= \sum_{\sigma\in \sum_n } {\rm sgn}(\sigma) \,\bl
z^{\bl m_\sigma},
$$
where $\bl z^{\bl m_\sigma}= z_1^{m_{\sigma(1)}}\cdots z_n^{m_{\sigma(n)}}$.
Thus for any $\bl p \in [\![n]\!]$, we have $$a_{\boldsymbol p}(\bl z) =
a_{\boldsymbol m + \boldsymbol \delta} (\bl z) =  \sum_{\sigma\in \sum_n } {\rm
sgn}(\sigma) \,\bl z^{(\boldsymbol m+ \boldsymbol \delta)_\sigma},$$
$\bl m \in [n]$ and it follows that
$$
a_{\bl p}(\bl z) = a_{\boldsymbol m+ \boldsymbol \delta}(\bl z) = \det\Big (
(\!( z_i^{p_j} )\!)_{i,j=1}^n \Big ),\, \bl p\in [\![n]\!].$$



The following Lemma clearly shows that the functions $a_{\bl p}$, $\bl  p\in
[\![n]\!]$, are orthogonal in the Hilbert space $\mathbb A^{(\lambda)}(\mathbb
D^n)$.
\begin{lem}\label{ortho}
The set
$S:=\{m_{\sigma(k)}-m_{\nu(k)}^\prime:\sigma,\nu\in\Sigma_n,m_i>m_j,
m_i^\i>m_j^\i \mbox{~for~}
i<j,m_1\neq m_1^\i,1\leq k\leq n\}$ $\not =$ $\{0\}.$
\end{lem}
\begin{proof} If there exist $\sigma,\nu\in \Sigma_n$ such that
$\sigma(k)=\nu(k)=1$ for some $k$, $1\leq k\leq n,$ then
$m_{\sigma(k)}-m_{\nu(k)}^\i=m_1-m_1^\i\neq 0.$
Therefore,  in this case, $S\not = \{0\}$.

Now suppose that there exists no $k, 1\leq k\leq n$, for which
$\sigma(k)=\nu(k)=1$. In this case, if possible, let $S=\{0\}.$
Fix $\sigma, \nu\in \Sigma_n$. Then there exists $k$ such that $\sigma(k)=1$ and
$\nu(k)=j>1.$ Now,
$m_{\sigma(k)}-m^\i_{\nu(k)}=m_1-m^\i_j.$ Pick $k^\i\neq k$ such that
$\sigma(k^\i)=j,\nu(k^\i)=\ell$, $\ell\neq j.$
Thus $m_{\sigma(k^\i)}-m^\i_{\nu(k^\i)}=m_j-m_\ell^\i.$ Choose $k^{\i\i}\neq k$
such that $\nu(k^{\i\i})=1, \sigma(k^{\i\i})=r>1$
and $m_{\sigma(k^{\i\i})}-m_{\nu(k^{\i\i})}=m_r-m^\i_1.$ However, we have
$m_1-m_j^\i=m_j-m_\ell^\i=m_r-m_1^\i=0.$ Clearly,
$m_r=m^\i_1>m_j^\i=m_1.$ Hence $m_r>m_1$ with $r>1,$ which is a contradiction.
\end{proof}

For $\lambda > 1$, the preceding Lemma says that the vectors $\bl z^{\bl
p_\sigma}$ are orthogonal, and hence the set $\{a_{\bl p}: \bl p \in
[\![n]\!]\}$
consists of mutually orthogonal vectors in $\mb A^{(\l)}(\D^n)$. The linear span
these vectors is dense in the Hilbert space
$\mb A_{\rm anti}^{(\l)}(\D^n).$ For $\bl p=(p_1,\ldots,p_n) \in [\![n]\!]$. The
norm of the vector
$a_{\bl p}$ is easily calculated:
\begin{eqnarray*}
\|a_{\bl p}\|_{\mathbb A^{(\lambda)}(\mathbb D^n)} &=& \Big \|\det\Big ( (\!(
z_i^{p_j} )\!)_{i,j=1}^n \Big ) \Big \|_{\mathbb A^{(\lambda)}(\mathbb D^n)}\\
&=&\Big \|\sum_{\sigma \in \Sigma_n} {\rm sgn}(\sigma) \prod_{k=1}^n
z_k^{p_{\sigma(k)}} \Big \|_{\mathbb A^{(\lambda)}(\mathbb D^n)}
=\sqrt\frac{n!\bl p!}{(\l)_{\bl p}},
\end{eqnarray*}
where
$\bl p!=\prod_{j=1}^n m_j!$ and $(\l)_{\bl p} =\prod_{j=1}^n(\l)_{m_j}.$
Putting $c_{\bl p}=\sqrt\frac{(\l)_{\bl p}}{{n!\bl p}!}$, we see that
$$\{e_{\bl p} = c_{\bl p} \,a_{\bl p} :\bl p \in [\![n]\!]\}$$ is an orthonormal
basis for $\mb A_{\rm anti}^{(\l)}(\D^n).$
So the reproducing kernel $K^{(\l)}_{\rm anti}$ for
$\mb A_{\rm anti}^{(\l)}(\D^n)$ is given by
$$K^{(\l)}_{\rm anti}(\bl z,
\bl w)=\d_{\bl p\in[\![n]\!]}e_{\bl p}(\bl z)\ov{e_{\bl p}(\bl
w)},\,\,\mbox{for}\,\, \bl z,\bl
w\in\D^n.$$
For all $\sigma\in\Sigma_n$, we have $e_{\sigma(\bl p)}(\bl z)\ov{e_{\sigma(\bl
p)}(\bl w)}=e_{\bl p}(\bl z)\ov{e_{\bl p}(\bl w)}$, $\bl z,\bl w\in\D^n$.
Therefore, it follows that
\bea\label{kernel} K^{(\l)}_{\rm anti}(\bl z, \bl w)=\d_{\bl
p\in[\![n]\!]}e_{\bl p}(\bl z)\ov{e_{\bl p}(\bl w)}=\frac{1}{n!}\d_{\bl p  \geq
0} e_{\bl p}(\bl z)\ov{e_{\bl p}(\bl w)},
\eea
where $\bl p \geq 0$ stands for all multi-indices $\bl p = (p_1,\ldots ,p_n) \in
\mathbb Z^n$ with the property that each $p_i \geq 0$ for $1\leq i \leq n$.

\begin{prop} The reproducing kernel $K^{(\l)}_{\rm anti}$ is given explicitly by
the formula:
$$K^{(\l)}_{\rm anti}(\bl
z, \bl w)=
\frac{1}{n!}
\det \Big (\big (\!\big ((1-z_j\bar w_k)^{-\l} \big )\!\big )_{j,k=1}^n\Big
),\,\, \bl z,\bl w\in\D^n.$$
\end{prop}
\begin{proof} 
For $\bl z, \bl w$ in $\mathbb D^n$, we have
\Bea \d_{\bl p \geq 0}e_{\bl p}(\bl z)\ov{e_{\bl p}(\bl w)}
&=& \frac{1}{n!}\d_{\bl p \geq 0}
\frac{(\l)_{\bl p}}{\bl p!} \,\det \Big ((\!(z_k^{p_j})\!)_{j,k=1}^n\Big )\det
\Big ((\!(\bar w_k^{p_j})\!)_{j,k=1}^n \Big )\\
&=&\frac{1}{n!}\d_{\bl p \geq 0}
\frac{(\l)_{\bl p}}{\bl p!} \,
\Big(\d_{\sigma\in\Sigma_n}\mbox{sgn}(\sigma)\displaystyle\prod_{i=1}^nz_i^{p_{
\sigma(i)}}\Big)
\Big (\d_{\nu\in\Sigma_n}\mbox{sgn}(\nu)\displaystyle\prod_{i=1}^n\bar
w_{\nu(i)}^{p_i}\Big)\\
&=& \frac{1}{n!}\d_{\bl p \geq 0}
\frac{(\l)_{\bl p}}{\bl p!}
\,\d_{\sigma,\nu\in\Sigma_n}\mbox{sgn}(\sigma)\mbox{sgn}(\nu)\displaystyle\prod_
{i=1}^n(z_i\bar
w_{\nu\sigma(i)})^{p_{\sigma(i)}}\\
&=&\frac{1}{n!}\d_{\sigma,\nu\in\Sigma_n}\mbox{sgn}(\nu \sigma)\d_{{\bl p} \geq
0}
\frac{(\l)_{\bl p}}{\bl p!} \,\displaystyle\prod_{i=1}^n(z_i\bar
w_{\nu\sigma(i)})^{p_{\sigma(i)}}\\
&=&\frac{1}{n!}\d_{\sigma,\nu\in\Sigma_n}\mbox{sgn}(\nu\sigma)\displaystyle
\prod_{i=1}^n(1-z_i\bar w_{\nu \sigma(i)})^{-\l}\\
&=&\frac{1}{n!}\d_{\psi\in\Sigma_n}\mbox{sgn}(\psi)\d_{\stackrel{\nu\sigma=\psi}
{\sigma,\nu\in\Sigma_n}
}\displaystyle\prod_{i=1}^n(1-z_i\bar
w_{\nu\sigma(i)})^{-\l}\\
&=& \d_{\psi\in\Sigma_n}\mbox{sgn}(\psi)\displaystyle\prod_{i=1}^n(1-z_i\bar
w_{\psi(i)})^{-\l}\\
&=& \det \Big ((\!((1-z_j\bar w_k)^{-\l})\!)_{j,k=1}^n\Big )\Eea
The desired equality follows from \eqref{kernel}.
\end{proof}

\subsection{Schur function}
The determinant function $a_{\bl m + \bl \delta}$ is divisible by each of the
difference $z_i-z_j$, $1\leq i < j \leq n$ and hence by the product
$$
\prod_{1\leq i < j \leq n} (z_i-z_j) =  \det \Big ( (\!( z_i^{n-j}
)\!)_{i,j=1}^n \Big )= a_{\boldsymbol \delta}(\bl z).
$$
The quotient $S_{\bl p}:=a_{\boldsymbol m+\boldsymbol \delta}/a_{\boldsymbol
\delta}$, $\bl p = \bl m + \bl \delta$, is therefore
well-defined and is called the Schur function \cite[pp.  40]{IGM}.
The Schur function $S_{\bl p}$ is symmetric and defines a function
on the symmetrized polydisc $\mathbb G_n$.
Since the Jacobian of the map $\mathbf s:\mathbb D^n \to \mathbb G_n$ coincides
with
$a_{\bl \delta}$, it follows from Lemma \ref{ortho} that the Schur functions
$\{S_{\bl p}:= a_{\bl m+\bl \delta}/a_{\bl \delta} :\bl p\in [\![n]\!]\}$
is a set of mutually orthogonal vectors in $\mb A^{(\l)}(\mathbb G_n)$.  The
linear span of these vectors
is dense in $\mb A^{(\l)}(\mathbb G_n).$ Also, the norms of these
vectors coincide with those of $a_{\bl p}$ in  $\mb A^{(\l)}(\mathbb G_n)$,
modulo the normalizing constant $\|J_\mathbf s\|_\lambda$, via the unitary map
$\Gamma$.  Hence $\| S_{\bl p}\|= \sqrt\frac{n!\bl p!}{\|J_\mathbf
s\|_\lambda(\l)_{\bl p }},$ $\bl p\in [\![n]\!].$
The set $\{\hat{e}_{\bl p}=c_{\bl p} \,S_{\bl p}:\bl p\in [\![n]\!]\}$ is an
orthonormal basis for $\mb A^{(\l)}(\mathbb G_n),$
where $c_{\bl p}=\sqrt\frac{\|J_\mathbf s\|_\lambda(\l)_{\bl p}}{n!\bl p!}.$
Thus we have proved:
\begin{thm}
For $\l > 0$, the reproducing kernel $\mathbf B_{\mb G_n}^{(\l)}$ for the
weighted Bergman space $\mb A^{(\l)}(\mathbb G_n)$ on the symmetrized poly-disc
is given by the formula:
\bea\label{skernel} \mathbf B^{(\l)}_{\mathbb G_n}(\mathbf s(\bl z), \mathbf
s(\bl w))
&=& \d_{\bl p \in [\![n]\!]}c^2_{\bl p}\, S_{\bl p}(\bl z)
\ov{S_{\bl p}(\bl w)}\\
&=&\frac{\|J_\mathbf s\|_\lambda^2}{n!}\frac{\det \Big( (\!((1-z_j\bar
w_k)^{-\l})\!)_{j,k=1}^n\Big)}{ a_{\boldsymbol{\delta}}(\bl z)
\overline{a_{\boldsymbol{\delta}}(\bl w)} }
\eea
for $\bl z,\bl w$ in $\mathbb D^n$.
\end{thm}
The case $\lambda = 2$ corresponds to the Bergman space on the symmetrized
polydisc. In this case, $\|J_\mathbf s\|_2=1$ and the formula for the the
Bergman kernel, except for the constant factor $\frac{1}{n!},$ was found in
\cite{EZ}.
(The factor $\frac{1}{n!}$ appears in our formula because we have chosen
the normalization $\|1\| =1$ for the constant function $1$ in the Hilbert space
$\mathbb A^{(\lambda)}(\mathbb G_n)$. However, as we will see below, it
disappears for the Hardy space on the symmetrized polydisc $\mathbb G_n$.)
However, the methods of this paper are very different form that of \cite{EZ},
and we hope it sheds some light on the nature of these kernel functions.
\begin{cor} The Bergman kernel on the symmetrized bidisc in $\mathbb C^2$ is
given by the formula
$$\mathbf B^{(2)}_{\mb G_2}(\bl
u, \bl v)=\frac{1}{2}\frac{2(1+u_2\bar v_2)-u_1\bar v_1}{((1-u_2\bar
v_2)^2-(u_1- u_2\bar v_1)(\bar v_1-\bar v_2 u_1))^2},$$ $\bl u=(u_1,u_2),\bl
v=(v_1,v_2)\in\mb G_2.$
\end{cor}
This corollary gives an explicit formula for the Bergman kernel function
for the symmetrized polydisc which is independent of the symmetriztion map
$\mathbf s$. It is possible to write down similar formulae for $n>2$ using
the Jacob-Trudy identity \cite[pp. 455]{FH}.

\section{The Hardy space and the Szeg\"{o} kernel for the symmetrized
polydisc}

Let $d\Theta$ be the normalized Lebesgue measure on the torus
$\mathbb T^n$, where $\mathbb T=\{\alpha :|\alpha| = 1\}$ is the unit circle.
Let $d\Theta_\mathbf s$ be the measure on the symmetrized polydisc $\mathbb G_n$
obtained by the change of variable formula:
$$
\int_{\partial \mathbb G_n} f \,d\Theta_\mathbf s
=\int_{\mathbb T^n} (f\circ\mathbf s)\, |J_\mathbf s|^2 d\Theta,
$$
where, as before, $J_\mathbf s(\bl z)$  is the complex Jacobian of the
symmetrization map $\mathbf s$.
The Hardy space $H^2(\mathbb G_n)$ on the symmetrized polydisc $\mathbb G_n$
consists of holomorphic functions on $\mathbb G_n$ with the property:
$$\mbox{sup}_{\,0<r<1}\int_{\mb T^n} |f\circ \mathbf s (r\,e^{i\Theta})|^2
|J_\mathbf s(r\,e^{i\Theta})|^2 d\Theta<\infty,\,\,e^{i\Theta} \in \mathbb
T^n.$$
We set the norm of $f\in H^2(\mb G_n)$ to be
$$\|f\|= \|J_{\bl s}\|^{-1}\Big \{\mbox{sup}_{0<r<1}\int_{\mb
T^n}|f\circ \mathbf s (r\,e^{i\Theta})|^2 |J_\mathbf s(r\,e^{i\Theta})|^2
d\Theta \Big \}^{1/2},$$
where $\|J_{\bl s}\|^2={{\int_{\mathbb T^n}}} |J_\mathbf s|^2 d\Theta.$  This
ensures, as before, $\|1\| = 1$.
Let $H^2(\D^n)$ be the Hardy space on the polydisc $\mathbb D^n$.
The operator $\Gamma :H^2(\mathbb G_n)\lo H^2(\D^n)$ given by
$\Gamma(f)= \|J_{\bl s}\|^{-1} J_\mathbf s \,(f\circ\mathbf s)$ for $f\in
H^2(\mathbb G_n)$ is then easily seen to be an
isometry. The subspace of anti-symmetric functions $H^2_{\rm anti}(\D^n)$ in the
Hardy space $H^2(\mathbb D^n)$ coincides with the image of $H^2(\mathbb G_n)$
under the isometry $\Gamma$. Thus the operator
$\Gamma:H^2(\mathbb G_n)\lo H^2_{\rm anti}(\D^n)$ is onto and therefore
unitary.

The functions $a_{\bl p}$, $\bl p\in [\![n]\!]$ continue to be an orthogonal
spanning set for the subspace
$H^2_{\rm anti}(\mathbb D^n)$.  All of the vectors $a_{\bl p}$ have the same
norm, namely, $\sqrt{n!}$.
Consequently, the set of vectors
$\{e_{\bl p}(\bl z):= \frac{1}{\sqrt{n!}}a_{\bl p}(\bl z):\bl p \in [\![n]\!]
\}$
is an orthonormal basis for the subspace $H^2_{\rm anti}(\mathbb D^n)$ of the
Hardy space on the polydisc,
while the set $\{\hat{e}_{\bl p}:= \frac{\|J_{\bl s}\|}{\sqrt{n!}} S_{\bl p}:
\bl p\in [\![n]\!]\}$ forms an orthonormal basis
for the Hardy space $H^2(\mathbb G_n)$ of the symmetrized polydisc $\mathbb G_n$
via the unitary map $\Gamma$. However, $\|J_{\bl s}\|=\sqrt{n!}$ and
consequently, $\hat{e}_{\bl p} = S_{\bl p}$.
Thus computations similar to the case $\lambda > 1$ yields an explicit
formula for the reproducing kernel
$K^{(1)}_{\rm anti}(\bl z, \bl w)$ of the subspace $H^2_{\rm anti}(\mathbb
D^n)$. Indeed,
$$K^{(1)}_{\rm anti}(\bl z,\bl w)=\frac{1}{n!}\det \Big((\!((1-z_j\bar
w_k)^{-1})\!)_{j,k=1}^n\Big ).$$
This is the limiting case, as $\lambda \to 1$.

Let $\mathbb S_{\mathbb G_n}$ be the S\"{z}ego kernel for the symmetrized
polydisc $\mathbb G_n$. Clearly,
$$
\mathbb S_{\mathbb G_n}(\mathbf s(\bl z),\mathbf s(\bl w))=
\frac{\det \Big ((\!((1-z_j\bar w_k)^{-1})\!)_{j,k=1}^n\Big )}{J_{\mathbf s}(\bl
z)\ov{J_{\mathbf s}(\bl w)}},\,\, \bl z,\bl w\in\D^n.
$$
Now, using  the well-known identity  due to Cauchy \cite[(4.3) pp 63]{IGM}, we
have
$$
\mathbb S_{\mathbb G_n}(\mathbf s(\bl z), \mathbf s(\bl w))=  \d_{\bl p \in
[\![n]\!]} S_{\bl p}(\bl z)
\ov{S_{\bl p}(\bl w)} = \prod_{j,k=1}^n (1-z_j\bar{w}_k)^{-1},\,\, \bl z,\bl w
\in \mathbb D^n.
$$
Therefore, we have a formula for the S\"{z}ego kernel of the symmetrized
polydisc $\mathbb G_n$, which we separately record below.


\begin{thm}
The S\"{z}ego kernel $\mathbb S_{\mathbb G_n}$ of the symmetrized polydisc
$\mathbb G_n$ 
is given by the formula
$$\mathbb S_{\mathbb G_n}(\mathbf s(\bl z),\mathbf s(\bl
w))=\displaystyle\prod_{j,k=1}^n (1-z_j\bar w_k)^{-1},
\bl z,\bl w\in\D^n.$$
\end{thm}

\section{An alternative approach to the computation of the kernel function}
Recall that the weighted Bergman space $\mathbb A^{(\lambda)}(\mathbb D^n)$ on
the polydisc $\mathbb D^n$ is the $n$-fold
tensor product $\otimes_{i=1}^n \mathbb A^{(\lambda)}(\mathbb D)$ of the
weighted Bergman spaces $\mathbb A^{(\lambda)}(\mathbb D)$
on the unit disc $\mathbb D$. The equivalence class $\widehat{\Sigma}_n$ of
finite dimensional irreducible representations
of the permutation group $\Sigma_n$ on $n$ symbols is parametrized by the
partitions $\bl p\in [n]$. Let $(V_{\bl p}, \bl p)$
be a representation corresponding to the partition $\bl p$. Then we have the
decomposition
$$
\mathbb A^{(\lambda)}(\mathbb D^n) = \oplus_{\bl p \in [n]} \mathbb
A^{(\lambda)}(\mathbb D^n, \bl p),
$$
where
$$
\mathbb A^{(\lambda)}(\mathbb D^n, \bl p) = \big \{f\in \mathbb
A^{(\lambda)}(\mathbb D^n, V_{\bl p}): \tau(\bl s) f(\bl s^{-1} \cdot \bl z) =
f(\bl z),  \bl s \in \Sigma_n \big \}
$$
and $\mathbb A^{(\lambda)}(\mathbb D^n, \bl p)\cong \mathbb
A^{(\lambda)}(\mathbb D, V_{\bl p}^\prime) \otimes V_{\bl p}^\prime$.
The orthogonal projection $\mathbb P_{\bl p}: \mathbb A^{(\lambda)}(\mathbb
D^n) \to \mathbb A^{(\lambda)}(\mathbb D^n, \bl p)$ is given by the formula
$$
(\mathbb P_{\bl p} f)(\bl z) = \frac{\chi_{\bl p}(1)}{n!} \sum_\tau \chi(\tau)
f(\tau^{-1}\cdot \bl z),
$$
where the sum is over all $\tau$ in $\Sigma_n$ and $\chi_{\bl p}$ is the
character
corresponding to the representation $V_{\bl p}$.
Schur orthogonality relations ensure that $\mathbb P_{\bl p}^2 = \mathbb P_{\bl
p}$ and it follows that $\mathbb P_{\bl p}$
is a projection.  Let $V_{\rm sgn}$ be the sign representation of the
permutation group $\Sigma_n$ and $\mathbb P_{\rm sgn}$ be the corresponding
projection.
\begin{thm}
The reproducing kernel $K^{(\lambda)}_{\rm sgn}$ of the Hilbert space $\mathbb
A^{(\lambda)}(\mathbb D^n, {\rm sgn})$ is given by the formula
\begin{eqnarray*}
K^{(\lambda)}_{\rm sgn} (\bl z, \bl w) &=& \big ( \mathbb P_{\rm sgn}\otimes
\mathbb P_{\rm sgn}^*\big ) \Big (\prod_{i=1}^n(1-z_i\bar{w}_i)^{-\lambda} \Big
)\\
&=& \frac{a_{\boldsymbol \delta}(\bl z) \overline{a_{\boldsymbol \delta}(\bl
w)}}{n!} \sum_{\bl p\in [\![ n]\!]} \frac{(\lambda)_{\bl m + \boldsymbol
\delta}}{(\bl m + \boldsymbol \delta)!} S_{\bl p}(z) \overline{S_{\bl p}(w)},
\end{eqnarray*}
where $S_{\bl p}$ is the Schur function with $\bl p = \bl m + \boldsymbol
\delta$.
\end{thm}
\begin{proof}
Recall that  $K^{(\lambda)}(\bl z, \bl w)= \sum_{\bl m \geq 0}^\infty
\frac{(\lambda)_{\bl m}}{{\bl m}!} (\bl z\bar{\bl w})^{\bl m}$, $\lambda > 1$,
is the reproducing kernel of the weighted Bergman spaces $\mathbb
A^{(\lambda)}(\mathbb D^n)$.  Therefore, we have
\begin{eqnarray*}
\big (\mathbb P_{\rm sgn} \otimes I\big ) K^{(\lambda)}_{\bl w}(\bl z)&=&
\sum_{\bl m \geq 0}^\infty \frac{(\lambda)_{\bl m}}{\bl m!} \bar{\bl w}^{\bl
m}\mathbb P_{\rm sgn} \big (\bl z^{\bl m}\big ).
\end{eqnarray*}
However, $\mathbb P_{\rm sgn} \big (\bl z^{\bl m} ) = \frac{1}{n!} \det  \Big
((\!(z_i^{m_j})\!)\Big)$ which is zero unless $\bl m$ is in the orbit under
$\Sigma_n$ of the weight $\bl p$ in $[\![n]\!]$. So, we conclude that
\begin{eqnarray*}
\big (\mathbb P_{\rm sgn} \otimes \mathbb P^*_{\rm sgn} \big )K^{(\lambda)}(\bl
z, \bl w) &=& \sum_{\bl m \geq 0}^\infty \frac{(\lambda)_{\bl m}}{\bl m!}
\mathbb P_{\rm sgn} \big (\bl z^{\bl m}\big ) \mathbb P_{\rm sgn} \big (
\bar{\bl w}^{\bl m} \big )\\
&=& \sum_{\bl p\in [\![ n ]\!]} \gamma_{\bl p} \frac{(\lambda)_{\bl p}}{\bl p !}
\, a_{\bl p}(\bl z)\overline{a_{\bl p}(\bl w)}\\
&=& a_{\boldsymbol \delta}(\bl z) \overline{a_{\boldsymbol \delta}(\bl w)}
\sum_{\bl p\in [\![ n]\!]} \gamma_{\bl p} \frac{(\lambda)_{\bl p }}{(\bl p)!}
S_{\bl p}(\bl z) \overline{S_{\bl p}(\bl w)}.
\end{eqnarray*}
It is then easy to see that $\gamma_{\bl p} = \frac{1}{n!}$ completing the
proof.
\end{proof}
Clearly, the two kernel functions $K^{(\lambda)}_{\rm sgn}$ and
$K^{(\lambda)}_{\rm anti}$ are equal.
As before, the kernel function $K^{(\lambda)}_{\rm sgn}$, via the unitary map
$\Gamma$, gives a kernel function for the
weighted Bergman spaces $\mathbb A^{(\lambda)}(\mathbb G_n)$ on the symmetrized
polydisc $\mathbb G_n$.
Further more, if $\lambda =1$, then
\begin{eqnarray*}
\mathbb S_{\mathbb G_n}\big (\bl s(\bl z), \bl s(\bl w)\big ) &=&
\frac{n!}{a_{\boldsymbol \delta}(\bl z) \overline{a_{\boldsymbol \delta}
(\bl w)}} K^{(1)}_{\rm sgn}(\bl z, \bl w)\\
&=&
\sum_{\bl p \in [\![n]\!]} S_{\bl p}(\bl z) \overline{S_{\bl p}(\bl w)}\\
&=&
\prod_{i,j=1}^n(1-z_i\bar{w}_j)^{-1}, \,\, \bl z, \bl w\in \mathbb D^n,
\end{eqnarray*}
where the last equality is the formula \cite[(4.3),  pp. 63]{IGM}.

\bibliographystyle{amsplain}

\end{document}